\documentclass[letterpaper,12pt]{article}

\usepackage[letterpaper,top=3cm,bottom=3cm,left=3cm,right=3cm,marginparwidth=1.75cm]{geometry}

\usepackage[english]{babel}
\usepackage[utf8x]{inputenc}
\usepackage[T1]{fontenc}

\usepackage{amsmath}
\usepackage{amsthm}
\usepackage{amssymb}
\usepackage[all,cmtip]{xy}
\usepackage{graphicx}
\usepackage{theoremref}

\theoremstyle{plain}
\newtheorem{theorem}{Theorem}[section]
\newtheorem{lemma}[theorem]{Lemma}
\newtheorem{prop}[theorem]{Proposition}

\theoremstyle{definition}

\usepackage{theoremref}

\newcommand{\Z}{\mathbb{Z}}

\newcommand{\K}{\mathbb{K}}

\newcommand{\BAG}{B^A(G)}
\newcommand{\BAH}{B^A(H)}

\newcommand{\BARG}{B^A_R(G)}
\newcommand{\BAKG}{B^A_{\K}(G)}
\newcommand{\MAG}{\mathcal{M}^A_G}

\newcommand{\DAG}{\mathcal{D}^A_G}

\newcommand{\SG}{\mathcal{S}_G}

\newcommand{\SGperf}{\mathcal{P}_G}

\title{Indecomposable factors of fibered Burnside rings}

\author{
Benjam\'in Garc\'ia\\
Centro de Ciencias Matem\'aticas\\
Universidad Nacional Aut\'onoma de M\'exico\\
Morelia, Michoac\'an 58089\\
Mexico\\
\texttt{benjamingarcia@matmor.unam.mx}
\and 
Alberto G. Raggi-Cárdenas\\
Centro de Ciencias Matem\'aticas\\
Universidad Nacional Aut\'onoma de M\'exico\\
Morelia, Michoac\'an 58089\\
Mexico\\
\texttt{agraggi@gmail.com}
}

\date{\today}

\begin{document}

\maketitle

\begin{abstract}
In this paper, we describe the indecomposable factors of a fibered Burnside ring in terms of bases coming from the standard basis. We provide a further characterization of each factor as the solvable component of the fibered Burnside ring of a corresponding Weyl group. 
\end{abstract}



\section*{Introduction} 

Fibered Burnside rings were introduced by Dress in his article \textit{"The ring of monomial representations I. Structure Theory"} \cite{Dress}, as a generalization of the Burnside ring $B(G)$ and the ring of monomial representations $D(G)$, with the aim of providing a common framework for several induction theorems. In Dress's article, a commutative ring $\Omega(G,A)$ called the \textit{ring of monomial representations of $G$ with coefficients in $A$} is defined as the Grothendieck ring of the category of finite $A$-fibered $G$-sets, where $G$ is a finite group and $A$ is a fixed finite $G$-module. Later, the name \textit{$A$-fibered Burnside ring of $G$} and the notation $\BAG$ have been preferred over $\Omega(G,A)$, with $A$ being an abelian group, as in most interesting cases the action of $G$ on $A$ is trivial. The species, prime spectrum, and primitive idempotents of $\BAG$ had already been determined in \cite{Dress}, with further refinements for particular cases given in \cite{Barker}, \cite{Boltje} and \cite{FotKul}. There are also explicit formulas for the primitive idempotents of the scalar extension $B^A_{\K}(G)$ in terms of the standard basis in \cite{BolYil}, where $\K$ is a "big enough" field. The conductors of the primitive idempotents of $B^{C_n}_{\K}(G)$ are known for a wide range of cases, including the Burnside ring and the ring of monomial representations \cite{GarRag}.

In this note, we provide bases for the indecomposable factors of $\BARG$ for an integral domain $R$ of characteristic $0$ in which no prime divisor of $|G|$ is invertible. In Section 1, we present the necessary background on the structure of $\BAG$ and the construction of its prime spectrum. Section 2 is devoted to the study of the indecomposable factors, where we also show that every indecomposable factor is isomorphic to the solvable component of a fibered Burnside ring of a corresponding Weyl group, just as in \cite[\S 3]{KiLuRa}. 

\subsection*{Notation and general setting}

Throughout this note, $G$ denotes a finite group and $A$ is a finite abelian group. We write $O^s(G)$ for the minimum normal subgroup $T\unlhd G$ such that $G/T$ is solvable, and we say that $G$ is \textit{perfect} if $O^S(G)=G$; it is known that $O^s(O^s(G))=O^s(G)$, therefore $O^s(H)$ is a perfect subgroup for each $H\leq G$. The set of subgroups of $G$ is denoted by $\SG$, the group $G$ acts by conjugation on $\SG$ and $[\SG]$ denotes a set of representatives of the conjugacy classes. For $H,K\leq G$, we write $H=_GK$ (resp. $H\leq_G K$) if $H$ is conjugate (resp. subconjugate) to $K$. For $H\leq G$, its \textit{Weyl group} $N_G(H)/H$ is denoted by $W_G(H)$. For an integer $n\geq 1$, the subgroup $\text{Tor}_n(A)\leq A$ consists of the elements $a\in A$ such that $a^n=1$.

\section{Background material}

Formally, the \textit{$A$-fibered Burnside ring} $\BAG$ is defined as the Grothendieck ring of the category $_G\text{set}^A$ of $A$-fibered $G$-sets having finitely many orbits, but we will be working with a more common presentation instead. An \textit{$A$-subcharacter} or \textit{$A$-monomial pair} of $G$ is a pair $(H,\phi)$ where $H\leq G$ and $\phi\in \text{Hom}(H,A)$. The set of all $A$-subcharacters is indicated by $\MAG$; the group $G$ acts on $\MAG$ by conjugation, and we write $[H,\phi]_G$ for the orbit of a pair $(H,\phi)$ and $N_G(H,\phi)$ for its stabilizer. We present the $A$-fibered Burnside ring $\BAG$ as the free abelian group on the orbit set $G\backslash \MAG$, which we refer to as the \textit{standard basis} or \textit{monomial basis}, for the multiplication defined by
 $$[H,\phi]_G[K,\psi]_G=\sum_{g\in[H\backslash G/K]}[H\cap{^gK},\phi\cdot{^g\psi}]_G$$
for $(H,\phi),(K,\psi)\in \MAG$, where $\phi\cdot {^g\psi}=\phi|_{H\cap{^gK}}{^g\psi|_{H\cap{^gK}}}$, with the basic element $[G,1]_G$ acting as a multiplicative identity. 

The group $G$ acts trivially on $\text{Hom}(G,A)$, and the group ring $\Z\text{Hom}(G,A)$ can be seen as a subring of $\BAG$ if we identify a morphism $\phi:G\longrightarrow A$ with $[G,\phi]_G$ in $\BAG$, with a retraction map
$$\pi_G:\BAG \longrightarrow \Z\text{Hom}(G,A), [H,\phi]_G\mapsto \begin{cases}
    \phi & \text{if }H=G,\\
    0 &\text{else},
\end{cases}$$
which is a ring homomorphism. There is also an embedding of the Burnside ring $B(G)$ in $\BAG$ sending $[G/H]$ to $[H,1]_G$, with a retraction $[H,\phi]_G\mapsto [G/H]$, both ring homomorphisms. The collection of rings $\left\{\BAH\mid H\leq G\right\}$ together with induction, restriction and conjugation maps form a Green functor over $G$ (see \cite[\S 2]{GarRag}). In particular, the \textit{restriction from $G$ to $H$}, defined as
$$\text{res}^G_H:\BAG \longrightarrow \BAH, [K,\phi]_G\mapsto \sum_{g\in [H\backslash G/K]}[H\cap {^gK},{^g\phi}|_{H\cap {^gK}}]_H,$$
is a ring homomorphism for all $H\leq G$.

If $R$ is a commutative ring, we write $B^A_R(G)$ for the scalar extension $R\otimes \BAG$ regarded as an $R$-algebra. We say that a field $\K$ of characteristic $0$ is \textit{big enough} for $\BAG$ if it contains a primitive root $1$ of order $\text{exp}(\text{Tor}_{\text{exp}(G)}(A))$. For a big enough field $\K$, let $\DAG$ be the set of pairs $(H,\Phi)$ with $H\leq G$ and $\Phi\in \text{Hom}(H,A)^{\star}:=\text{Hom}(\text{Hom}(H,A),\K^{\times})$, on which $G$ by conjugation again; we write $[H,\Phi]_G$ for the orbit of $(H,\Phi)$ and $N_G(H,\Phi)$ for its stabilizer. For each $(H,\Phi)\in \DAG$, the composite map
$$\xymatrix{
s^{G,A}_{H,\Phi}:\BAG \ar[rr]^{\text{res}^G_H} &&\BAH \ar[rr]^-{\pi_H} &&\Z \text{Hom}(H,A) \ar[rr]^-{\Phi} &&\K
},$$
where $\Phi$ is extended linearly to $\Z\text{Hom}(H,A)$, is a ring homomorphism. Explicitly, for $(K,\psi)\in \MAG$, we have
$$s^{G,A}_{H,\Phi}\left([K,\psi]_G\right)=\sum_{\substack{g\in [G/K],\\
H\leq {^gK}}}\Phi\left({^g\psi|_H}\right);$$
in particular, $s^{G,A}_{H,\Phi}\left([K,\psi]_G\right)=0$ whenever $H$ is not subconjugate to $K$. By \cite[Theorem 3.1]{BolYil}, the product of the extended maps
$$\prod_{(H,\Phi)\in [\DAG]}s^{G,A}_{H,\Phi}:B^A_{\K}(G) \longrightarrow \prod_{(H,\Phi)\in [\DAG]}\K$$
is an isomorphism of $\K$-algebras, where $[\DAG]$ stands for a set of representatives of the conjugacy classes of $\DAG$. This implies that the maps $s^{G,A}_{H,\Phi}$ account for all ring homomorphisms from $\BAG$ to $\K$, whose images lie in $\Z[\zeta]$, where $\zeta$ is a root of $1$ of order $\text{exp}\left(\text{Tor}_{\text{exp}(G)}(A)\right)$. Then the primitive idempotents of $\BAKG$  are elements $e^{G,A}_{H,\Phi}$ for $[H,\Phi]_G\in G\backslash \DAG$, determined by the property
$$s^{G,A}_{K,\Psi}\left(e^{G,A}_{H,\Phi}\right)=\begin{cases}
    1 &\text{if }[K,\Psi]_G=[H,\Phi]_G,\\
    0 &\text{else,}
\end{cases}$$
for $(K,\Psi)\in \DAG$. By \cite[Theorem 3.2]{BolYil}, the idempotent $e^{G,A}_{H,\Phi}$ can be expressed as
\begin{align}
    e^{G,A}_{H,\Phi}=\frac{1}{|N_G(H,\Phi)||\text{Hom}(H,A)|}\sum_{\substack{
    K\leq H,\\
    \phi\in \text{Hom}(H,A)
    }
    }|K|\mu(K,H)\Phi(\phi^{-1})[K,\phi|_K]_G,
\end{align}
where $\mu$ is the Möbius function on the poset of subgroups of $G$.

We recall the construction of the prime spectrum of $\BARG$ when $R$ is an integral domain of  characteristic zero; this was done by Dress in \cite{Dress} in a quite general setting, but we follow the technique of Fotsing \& Kulshammer in \cite{FotKul}. Let $\zeta$ be a primitive $n$-th root of $1$ in some algebraic extension $\K$ of $Q(R)$, where $n=\text{exp}(\text{Tor}_{\text{exp}(G)}(A))$. For every $(H,\Phi)\in \DAG$, the extended map $s^{G,A}_{H,\Phi}:B^A_{R[\zeta]}(G)\longrightarrow R[\zeta]$ is surjective, and since
$$\bigcap_{(H,\Phi)\in [\DAG]} \text{ker}\left(s^{G,A}_{H,\Phi}\right)=0$$
is a finite intersection of prime ideals, we have that every prime ideal of $B^A_{R[\zeta]}(G)$ has the form
$$\mathbf{P}(H,\Phi,P):=\left(s^{G,A}_{H,\Phi}\right)^{-1}(P)$$
for some $(H,\Phi)\in \DAG$ and $P\in \text{Spec}(R[\zeta])$. Since $B^A_{R[\zeta]}(G)$ is integral over $\BARG$, we get a surjective application
$$\DAG \times \text{Spec}\left(R[\zeta]\right)\longrightarrow \text{Spec}\left(\BARG\right),$$
sending $(H,\Phi,P)$ to the prime ideal
$$\mathfrak{p}(H,\Phi,P):=\mathbf{P}(H,\Phi,P) \cap \BARG.$$
We will not discuss when $\mathfrak{p}(H,\Phi,P) = \mathfrak{p}(K,\Psi,Q)$, but a classical characterization can be achieved following the techniques in \cite{FotKul}. 

If no prime divisor of $|G|$ is a unit in $R$, then $\text{Spec}\left(\BARG\right)$ and $\text{Spec}\left(B_R(G)\right)$ have the same number of connected components. By \cite[Corollary 2]{Dress}, if $J$ is a perfect subgroup of $G$, the subset
$$\mathcal{K}_J=\left\{\mathfrak{p}(H,\Phi,P) \in \text{Spec}\left(\BARG\right) \mid O^s(H)=_GJ\right\}$$
is a connected component of $\text{Spec}\left(\BARG\right)$, and we have $\mathcal{K}_J\neq \mathcal{K}_I$ whenever $I$ is another perfect subgroup such that $I\neq_GJ$. Therefore, the connected components of $\text{Spec}\left(\BARG\right)$ are in bijection with the conjugacy classes of perfect subgroups of $G$, and their number does not depend on $A$. By \cite[Proposition 3.1]{Yoshida}, the primitive idempotents of $B_R(G)$ are the elements
$$e^G_{[J]}=\sum_{\substack{H\in [\SG],\\
O^s(H)=_G J}}e^G_H,$$
where $J$ runs over the perfect subgroups of $G$ up to conjugation. These idempotents remain primitive in $\BARG$ by a simple counting argument.

\section{Indecomposable factors}

Let $R$ be an integral domain of characteristic $0$ such that no prime divisor of $|G|$ is in $R^{\times}$. We take a big enough field extension $\K$ of the field of fractions of $R$, and consider the natural embedding $\BARG$ into $\BAKG$. We write $\SGperf$ for the set of perfect subgroups of $G$ and $[\SGperf]$ for a set of representatives of its conjugacy classes.

\begin{lemma}\thlabel{blocks}
    For $J\in [\SGperf]$, we have
    \begin{equation}
        e^G_{[J]}=\sum_{\substack{
    (H,\Phi)\in [\DAG],\\
    O^s(H)=_GJ
    }} e^{G,A}_{H,\Phi}.
    \end{equation}
\end{lemma}
\begin{proof}
    For $H\leq G$ and $(K,\Phi)\in \DAG$, we have $s^{G,A}_{K,\Phi}\left([H,1]_G\right)=\varphi_K\left(G/H\right)$ (the $K$-th mark on $G/H$). Because  $e^G_{[J]}$ can be expressed as a sum of terms of the form $[H,1]_G$, we have
    $$s^{G,A}_{K,\Phi}\left(e^G_{[J]}\right)=\varphi_K\left(e^G_{[J]}\right)=\begin{cases}
        1 &\text{if }O^s(K)=_GJ,\\
        0 &\text{else}
    \end{cases}$$
    for each $(K,\Phi)\in \DAG$.
\end{proof}

From the previous discussion on primitive idempotents, we get a decomposition
\begin{align}
    \BARG=\bigoplus_{J\in [\SGperf]}\BARG e^G_{[J]}
\end{align}
as a sum of indecomposable factors or components. The factor corresponding to $J=1$ is known as the \textit{solvable component} of $\BARG$.

\begin{prop}\thlabel{Prop1}
    For $J\in [\SGperf]$, we have
    \begin{align}
        \BARG e^{G}_{[J]}=\bigoplus_{
        \substack{(H,\phi)\in [\MAG],\\
        O^s(H)=_GJ}}R [H,\phi]_G e^G_{[J]}.
    \end{align}
    In particular,
    \begin{align}
        \BARG e^{G}_{[1]}=\bigoplus_{
        \substack{(H,\phi)\in [\MAG],\\
        O^s(H)=1}}R [H,\phi]_G .
    \end{align}
\end{prop}
\begin{proof}
    It suffices to prove it for $R=\Z$. We first show that
    \begin{align}
        \BAG e^G_{[J]}= \sum_{\substack{
    (H,\phi)\in [\MAG],\\
    O^s(H)=_GJ
    }} \Z [H,\phi]_Ge^G_{[J]}.
    \end{align}
    The right-hand side is clearly contained on the left-hand side. From Equation (1) and since ${e^G_{[J]}}^2=e^G_{[J]}$, we have
    $$e^G_{[J]}=\sum_{\substack{
    (H,\phi)\in[\MAG],\\
    O^s(H)\leq_G J}}a_{H,\phi}[H,\phi]_G=\sum_{\substack{
    (H,\phi)\in[\MAG],\\
    O^s(H)\leq_G J}}a_{H,\phi}[H,\phi]_Ge^G_{[J]}$$
    for some $a_{H,\phi}\in\Z$ after grouping similar terms. From the multiplication formula, for each $(H,\phi)\in \MAG$ we have
    \begin{align}
        [H,\phi]_Ge^G_{[J]}=\sum_{\substack{
    (K,\psi)\in[\MAG],\\
    O^s(K)\leq_G J}}b_{K,\psi}[K,\psi]_Ge^G_{[J]}
    \end{align}
    for some $b_{K,\psi}\in \Z$ after grouping similar terms again, showing that it suffices to consider only $[H,\phi]e^G_{[J]}$ with $O^s(H)\leq_GJ$. If $O^s(H)$ is a proper subconjugate of $J$, we consider two situations for arbitrary $(L,\Psi)\in\DAG$: if $O^s(L)\neq_G J$ then $s^{G,A}_{L,\Psi}\left(e^G_{[J]}\right)=0$; if $O^s(L)=_GJ$, then $L$ cannot be subconjugate to $H$ and so $s^{G,A}_{L,\Psi}\left([H,\phi]_G\right)=0$. Hence, $s^{G,A}_{L,\Psi}\left([H,\phi]_Ge^G_{[J]}\right)=0$ for all $(L,\Psi)$ whenever $O^s(H)\lneq_GJ$, thus $[H,\phi]_Ge^G_{[J]}=0$. This shows that we can take $b_{K,\psi}=0$ in (7) whenever $O^s(K)\neq_GJ$, proving Equation (6). To see that the sum is direct, consider an expression
    $$0=\sum_{\substack{
    (H,\phi)\in [\MAG],\\
    O^s(H)=_GJ
    }}b_{H,\phi}[H,\phi]_Ge^G_{[J]}$$
    for some $b_{H,\phi}\in \Z$. Suppose that some $b_{H,\phi}$ are nonzero and take $(K,\psi)$ with $K$ of maximal order such that $b_{K,\psi}\neq 0$. For all $\Phi\in \text{Hom}(K,A)^{\star}$, we have
    \begin{align*}
        0 &=s^{G,A}_{K,\Phi}\left(\sum_{\substack{
    (H,\phi)\in [\MAG],\\
    O^s(H)=_GJ
    }}b_{H,\phi}[H,\phi]_Ge^G_{[J]}\right)\\
    &=\sum_{\substack{
    (H,\phi)\in [\MAG],\\
    H=_GK
    }}b_{H,\phi}s^{G,A}_{K,\Phi}\left([H,\phi]_G\right)\\
    &=\sum_{\phi\in [N_G(K)\backslash \text{Hom}(K,A)]}b_{K,\phi}\sum_{gK\in N_G(K)/K}\Phi\left({^g\phi}\right)\\
    &=\Phi\left(\sum_{\phi\in \text{Hom}(K,A)}b_{K,\phi}[N_G(K,\phi):K]\phi\right),
    \end{align*}
    hence 
    $$\sum_{\phi\in \text{Hom}(K,A)}b_{K,\phi}[N_G(K,\phi):K]\phi=0$$ 
    in $\Z\text{Hom}(K,A)$, and so $b_{K,\phi}[N_G(K,\phi):K]=0$ for all $\phi$; in particular $b_{K,\psi}=0$, a contradiction. Finally, for $J=1$ and $(H,\phi)$ with $H$ solvable, we have 
    \begin{align*}
        s^{G,A}_{K,\Phi}\left([H,\phi]_Ge^{G}_{[1]}\right)= s^{G,A}_{K,\Phi}\left([H,\phi]_G\right)
    \end{align*}
    for all $(K,\Phi)\in \DAG$, proving $[H,\phi]_Ge^{G}_{[1]}=[H,\phi]_G$.
\end{proof}

\begin{prop}\thlabel{Prop2}
    For $J\in [\SGperf]$, we have
    \begin{align}
        \BARG e^{G}_{[J]}\cong B^A_R(W_G(J))e^{W_G(J)}_{[1]}.
    \end{align}
\end{prop}
\begin{proof}
    We prove it for $R=\Z$. First, note that if $J\leq H\leq N_G(J)$, then $O^s(H/J)=J/J$ if and only if $O^s(H)=J$. We define a map
    $$f:\bigoplus_{\substack{
    (H/J,\phi)\in [\mathcal{M}^A_{W_G(J)}],\\
    O^s(H)=J}}
    \Z[H/J,\phi]_{W_G(J)}
    \longrightarrow \bigoplus_{\substack{
    (H,\tilde{\phi})\in [\MAG],\\
    O^s(H)=_GJ}} \Z [H,\tilde{\phi}]_Ge^G_{[J]}
    $$
    sending $[H/J,\phi]_{W_G(J)}$ to $[H,\tilde{\phi}]_G$, where $\tilde{\phi}$ is the inflation from $H/J$ to $H$ of $\phi$. This is easily verified to be well-defined. It is onto since each pair $(H,\phi)$ with $O^s(H)=_GJ$ is conjugate to some $(K,\psi)$ with $O^s(K)=J$; for this pair, we have $K\leq N_G(J)$ and $\psi$ is the inflation of some $\widehat{\psi}:K/J\longrightarrow A$ since $J\leq K'\leq \text{ker}(\psi)$, therefore $f\left([K/J,\widehat{\psi}]_{W_G(J)}\right)=[H,\phi]_G$. We check injectivity on the basis: if $[H,\widetilde{\phi}]_Ge^G_{[J]}=[K,\widetilde{\psi}]_Ge^G_{[J]}$ for some $(H/J,\phi),(K/J,\psi)\in \mathcal{M}^A_{W_G(J)}$, then $[H,\widetilde{\phi}]_G=[K,\widetilde{\psi}]_G$ by \thref{Prop1}, hence there is some $g\in G$ such that $^g(H,\widetilde{\phi})=(K,\widetilde{\psi})$. This implies that $J=O^s(K)=O^{s}({^gH})={^gJ}$, and so $g\in N_G(J)$, hence ${^{\overline{g}}(H/J,\phi)}=(K/J,\psi)$ for $\overline{g}=gJ \in W_G(J)$. This shows that $f$ is an isomorphism of abelian groups. We now prove that $f$ is multiplicative: for $(H/J,\phi),(K/J,\psi)\in \mathcal{M}^A_{W_G(J)}$, we have
    \begin{align*}
        f\left([H/J,\phi]_{W_G(J)}[K/J,\psi]_{W_G(J)}\right) &=\sum_{\overline{g}\in [(H/J)\backslash W_G(J)/(K/J)]}f\left([(H/J)\cap {^{\overline{g}}(K/J)},\phi\cdot {^{\overline{g}}\psi}]_{W_G(J)}\right)\\
        &=\sum_{g\in [H\backslash N_G(J)/K]}[H\cap {^gK},\widetilde{\phi}\cdot {^g\widetilde{\psi}}]_Ge^G_{[J]}.
    \end{align*}
    To see that the latter sum is the product of the images, note that if $g\in G$ is such that $[H\cap {^gK},\widetilde{\phi}\cdot {^g\widetilde{\psi}}]_Ge^G_{[J]}\neq 0$, then $O^s(H\cap {^gK})=_GJ$. But $O^s(H\cap{^gK})\leq J\cap{^gJ}$, which implies $J={^gJ}$ and $g\in N_G(J)$. Hence
    \begin{align*}
        \sum_{g\in [H\backslash N_G(J)/K]}[H\cap {^gK},\widetilde{\phi}\cdot {^g\widetilde{\psi}}]_Ge^G_{[J]} &=\sum_{g\in [H\backslash G/K]}[H\cap {^gK},\widetilde{\phi}\cdot {^g\widetilde{\psi}}]_Ge^G_{[J]}\\
        &=\left([H,\widetilde{\phi}]_Ge^G_{[J]}\right)\left([K,\widetilde{\psi}]_Ge^G_{[J]}\right)\\
        &=f\left([H/J,\phi]_{W_G(J)})f([K/J,\psi]_{W_G(J)}\right).
    \end{align*}
    This also shows that $f$ sends $e^{W_G(J)}_{[1]}$ to $e^G_{[J]}$ since $f\left(e^{W_G(J)}_{[1]}\right)$ is a nonzero idempotent in $\BAG e^G_{[J]}$, which is indecomposable. Hence $f$ is a ring isomorphism.
\end{proof}

It is known that a group $G$ is solvable if and only if $\BAG$ is indecomposable, and in the case of odd order, proving the Feit-Thompson Theorem is equivalent to proving that $[G,1]_G$ is a primitive idempotent in $\BAG$, which exhibits the difficulty of the second assertion rather than providing an easier alternative proof. From \cite[Theorem 5.3]{KiLuRa}, $B(G)$ determines $G$ (up to isomorphism) for a wide range of simple groups; the proof relies strongly on \cite[Theorem 3.1]{KiLuRa} (that we have generalized here) and the existence of normalizing isomorphisms. However, the existence of normalizing isomophisms for fibered Burnside rings is an open problem.

\subsection*{Acknowledgments}

The first author was supported by SECIHTI with an "Estancias Posdoctorales por México" grant.

\end{document}